\documentclass[draft]{article}
\title{The Noncommutative Geometry Generalization of Fundamental Group}
\author{Petr R. Ivankov and Nickolay P. Ivankov }
\date{April 24, 2006}

\usepackage{amsmath,amsthm, pb-diagram, lamsarrow, pb-lams}

\chardef\bslash=`\\ 

\newtheorem{thm}{Theorem}[section]

\newtheorem{lem}[thm]{Lemma}

\newtheorem{ax}{Axiom}

\theoremstyle{definition}
\newtheorem{defn}{Definition}[section]

\theoremstyle{remark}

\newcommand{\A}{\mathcal{A}}
\newcommand{\B}{\mathcal{B}}

\newcommand{\wt}{\widetilde}
\newcommand{\wh}{\widehat}

\newcommand{\eval}[2][\right]{\relax
  \ifx#1\right\relax \left.\fi#2#1\rvert}

\begin{document}
\maketitle
\markboth{Petr R. Ivankov, Nickolay P. Ivankov}{The Noncommutative Geometry Generalization of Fundamental Group}

\renewcommand{\sectionmark}[1]{}

\begin{abstract}
{A notion of fundamental group of spectral triples \cite{varilly:noncom} has been introduced. The notion uses a noncommutative analogue of
unramified coverings. It was shown that in commutative case this fundamental group is a profinite completion of fundamental group of corresponding
Riemann manifold.  }  
\end{abstract}

\section{Introduction}

In algebraic topology the fundamental group may be defined using closed paths or automorphisms of universal covering.
However, noncommutative geometry spectral triple has no paths and even points in general case. 
Similary in algebraic geometry we can not always find a set of paths that provides a good definition of fundamental group.
But it is possible to define unramified coverings. These coverings enable us to define analogue of fundamental group in 
algeraic geometry.	The algebraic geometry fundamental group of complex algebraic manifold is the profinite completion of algebraic 
topology fundamental group \cite{milne:etale}.
The notion of ramified covering in noncommutative geometry was introduced in \cite{connesdebois:3dsphere}.
Using very similar approach we have introduced the notion of unramified finite covering in noncommutative geometry. The spectral triple
coverings enable us to define the noncommutative version of fundamental group. If a spectral triple is 
commutative then its fundamental group is the profinite completion of the fundamental group of Riemann manifold that is associated to the triple.

\section{Preliminaries}
Recall definition of real spectral triple\cite{varilly:noncom}. A real spectral triple consits
of a set of four (five) objects $(\A, \mathcal{H}, D, J (,\Gamma))$, of following types:

 (1) $\A$ is a $pre-C^* algebra$;
 
 (2) $\mathcal{H}$ is a $Hilbert$ $space$ carrying faithful representation $\pi$ of algebra $\A$;
 
 (3) $D$ is a $selfadjoint$ $operator$ on $\mathcal{H}$ with compact resolvent; 

 (4) $J$ is a $antilinear$ $isometry$ of $\mathcal{H}$ onto itself; 

 (5) $\Gamma$ is a $selfadjoint$ $unitary$ $operator$ on $\mathcal{H}$ so that $\Gamma^2=1$.
 
 If  $\Gamma$ is present we say that triple is $even$, otherwise it is $odd$.
The reader can find complete description of spectral triples, and those axioms at \cite{varilly:noncom}.
We will use following notation of spectral triples $\mathbf{A}=(\A, \mathcal{H}, D, J ,\Gamma)$
or  $\mathbf{A}=(\A, \mathcal{H}, D, J)$ or $\mathbf{A}=(\A, \mathcal{H}, D)$.

In the following text $D$, $J$ ,$(\Gamma)$ will be called as ``operators".

Recall the notion of unitary equivalence of geometries \cite{varilly:noncom}. Unitary equivalence of spectral triple
 $\mathbf{A}=(\A, \mathcal{H}, D)$ is associated with such unitary operator $U$ on $\mathcal{H}$
that $U\pi(\A)U^{-1}=\pi(\A)$.  

Let $G(\mathbf{A})$ or $G(\A, \mathcal{H}, D)$ be the group of unitary equivalences of the triple $\mathbf{A}=(\A, \mathcal{H}, D)$.
Every unitary equivalence $U$ defines * - automoprphism $\sigma$ of $\A$ that satisfies following condition $U\pi (a)U^{-1}=\pi (\sigma (a))$.
Hence $G(\mathbf{A})$ acts on $\mathcal{H}_A$ and $\A$.

The orientability axiom \cite{varilly:noncom} assume the existence of fundamental Hochschild cycle $c_A\in Z_n(\A\otimes\A^0)$. 
Every homomorhism $f:\A\rightarrow\B$ of algebras defines
natural transformation $f^*$ of Hochschild cycles. Thus we have $f^*(c_A)\in Z_n(\B\otimes\B^0)$.

\section{Finite coverings. Fundamental group}
In this section we shall define noncommutative analogue of covevings using sructures of spectral triple only. We will
show below that in commutative case there is one to one correspondence between these coverings and finitely sheeted coverings of corresponding
Riemann manifold. These coverings enable to define fundamental group.
\subsection{Elements of finite covering}

Let $G \subset G(\mathbf{A})$ be a finite subgroup of $G(\mathbf{A})$. Elements of $G$ act on $\A$ and $\mathcal{H}$.
Thefore group $G$ defines projector $P$ on $\A$ 
and $\mathcal{H}$ by following formula:
\begin{equation}\label{surgroup}
P_G=\frac{1}{|G|}\sum_{g\in G} g
\end{equation}
The image of $P$ is a subalgebra (subspace) $\A^G$ ($\mathcal{H}^G$) of $\A$ ($\mathcal{H}$).

\begin{defn}\label{defn-per}
Finite covering of spectral triples $p:(\A, \mathcal{H}_A, D_A)\rightarrow (\B, \mathcal{H}_B, D_B)$ consists of such pair $(p_1,p_2)$
of injective *-homomorphism $p_1:\B\rightarrow\A$ and homomorphism $p_2:\mathcal{H}_B\rightarrow\mathcal{H}_A$ that:

(i) There exists a finite subgroup $G(\mathbf{A}, \mathbf{B}) \subset G(\A, \mathcal{H}, D)$ that 
image of $p_1$ ($p_2$) is $\A^{G(\mathbf{A}, \mathbf{B})}$ ($\mathcal{H}^{G(\mathbf{A}, \mathbf{B})}$);

(ii) Homomorpisms $p_1$ and $p_2$ satisfy to the condition 

$p_2(\pi(b)h)=\pi(p_1(b))p_2(h)$ for
any $b\in\B$ and $h\in\mathcal{H}_B$;

(iii) Spectral triples $\mathbf{A} = (\A, \mathcal{H}_A, D_A)$ , $\mathbf{B} = (\B, \mathcal{H}_B, D_B)$
and group $G(\mathbf{A}, \mathbf{B})$ satisfy following axioms 1 - 7.

\end{defn}

\subsection{Axioms of finite covering}
\begin{ax}
$\A$ is a finitely generated projective $\B-$ module.
\end{ax}
\begin{ax}
There is the natural isomorphism of $\A-$ modules $\mathcal{H}_A \approx \A\otimes_{\B} \mathcal{H}_B$
\end{ax}
\begin{ax}
There exists such surjective homomorphism $\phi_{AB}:G(\mathbf{A})\rightarrow G(\mathbf{B})$ that the following diagram is commutative
\[
\begin{diagram}
	\node{\mathcal{H}_A} \arrow[2]{e,t}{P_{G(\mathbf{A}, \mathbf{B})}} \arrow{s,l}{g}
								\node[2]{\mathcal{H}_B} \arrow{s,r}{\phi_{AB}(g)} \\
	\node{\mathcal{H}_A} \arrow[2]{e,b}{P_{G(\mathbf{A}, \mathbf{B})}}	\node[2]{\mathcal{H}_B}
\end{diagram}
\]
where $P$ is the projection defined by equation \eqref{surgroup}.
\end{ax}
\begin{ax}
Every element $g\in\ G(\mathbf{A})$ that is identity on  $p_2(\mathcal{H}_B)$ belongs to $G(\mathbf{A}, \mathbf{B})$.
\end{ax}
\begin{ax}
Every $g\in G(\mathbf{A}, \mathbf{B})$ commute with $D_A$, $J_A$ (, $\Gamma_A$)
and $D_B$, $J_B$ (, $\Gamma_B$) are restrictions of $D_A$, $J_A$ (, $\Gamma_A$) on $\mathcal{H}_B$.
\end{ax}
\begin{ax}
The dimension of the triple $\mathbf{A} = (\A, \mathcal{H}_A, D_A)$ equals to the dimension of the triple $\mathbf{B} = (\B, \mathcal{H}_B, D_B)$.
Definition of a triple dimension is contained in \cite{varilly:noncom}.
\end{ax}
\begin{ax}
If $c_A\in Z_n(\A\otimes\A^0)$ is fundamental cycle on $\mathbf{A}$ and $c_B\in Z_n(\B\otimes\B^0)$ is fundamental cycle on $\mathbf{B}$ then
$c_A=f^*(c_B)$.  
\end{ax}

In the following text we consider $\B$ as subalgebra of $\A$ and $\mathcal{H}_B$ as subspace of 
$\mathcal{H}_A$.

\subsection{Composition of finite coverings}
\begin{lem}\label{lem-comp}
Let
\[
\begin{diagram}
	\node{\mathbf{A}= (\A, \mathcal{H}_A, D_A)} \arrow {e,t}{f} \node {\mathbf{B} = (\B, \mathcal{H}_B, D_B)} \arrow{e,l}{g} \node {\mathbf{C}= (\mathcal{C}, \mathcal{H}_C, D_C)}
\end{diagram}
\]
be a diagram with finite coverings of spectral triples, and $f=(f_1, f_2)$, $g=(g_1, g_2)$. Then the pair $(f_1g_1, f_2g_2)$ defines finite covering
$gf:\mathbf{A}\rightarrow\mathbf{C}$

\end{lem}
\begin{proof}
Let us check (i) - (iii)
(i) Let $ G(\mathbf{A}, \mathbf{C}) \subset G(\mathbf{A})$ be the subroup of those automorphisms that are identities on $p_1(\mathcal{C})$ and $p_2(\mathcal{H}_C)$ 
According to axiom 4 we have a surjective homomorphism $G(\mathbf{A}, \mathbf{C}) \rightarrow G(\mathbf{B}, \mathbf{C})$
and even an exact sequence 

\begin{equation}\label{exact-covering}
\{e\} \rightarrow G(\mathbf{A}, \mathbf{B}) \rightarrow G(\mathbf{A}, \mathbf{C}) \rightarrow G(\mathbf{B}, \mathbf{C})\rightarrow \{e\}
\end{equation}

Since $G(\mathbf{A}, \mathbf{B})$ and $G(\mathbf{B}, \mathbf{C})$ are finite then $G(\mathbf{A}, \mathbf{C})$ is finite.
(ii) Follows from simple direct calculation and omitted here.
(iii) 

Axiom 1. The $\A$ is the direct summand of $\mathcal{B}^m$, and the $\B$ is direct summand of $\mathcal{C}^n$.
Hence $\A$ is the direct summand of $\mathcal{C}^{mn}$.

Axiom 2. This axiom follows from next equivalences:

$\mathcal{H}_B \approx \B\otimes_{\mathcal{C}} \mathcal{H}_C$

$\mathcal{H}_A \approx \A\otimes_{\B}\mathcal{H}_B \approx\A\otimes_{\mathcal{B}}\B\otimes_{\mathcal{C}}\mathcal{H}_C\approx \A\otimes_{\mathcal{C}} \mathcal{H}_C$

Axiom 3.
Follows from the next commutative diagram:

\[
\begin{diagram}
	\node{\mathcal{H}_A}\arrow{s, l}{g} \arrow[2]{e,t}{P_{G(\mathbf{A}, \mathbf{B})}} 
								\node[2]{\mathcal{H}_B} \arrow{s,r}{\phi_{AB}(g)}\arrow[2]{e,t}{P_{G(\mathbf{B}, \mathbf{C})}} 
								\node[2]{\mathcal{H}_C} \arrow{s,r}{\phi_{AC}(g)} \\
	\node{\mathcal{H}_A} \arrow[2]{e,b}{P_{G(\mathbf{A}, \mathbf{B})}}	\node[2]{\mathcal{H}_B}
	\arrow[2]{e,b}{P_{G(\mathbf{B}, \mathbf{C})}}	\node[2]{\mathcal{H}_C}
\end{diagram}.
\]

Axiom 4.
Follows from the exact sequence \eqref{exact-covering} of finite groups.

Axiom 5. 
Let $g\in G(\mathbf{A}, \mathbf{C})$,
$\mathcal{H}_1$ is ortogonal supplement of $\mathcal{H}_C$ in $\mathcal{H}_B$,
and $\mathcal{H}_2$ is ortogonal supplement of $\mathcal{H}_B$ in $\mathcal{H}_A$.
Then $\mathcal{H}_A=\mathcal{H}_C\oplus\mathcal{H}_1\oplus\mathcal{H}_2$
and $\mathcal{H}_B=\mathcal{H}_C\oplus\mathcal{H}_1$

The $\phi_{AB}(g)$ is represented by operator on $\mathcal{H}_B$ that looks like:

\begin{equation}\label{diag1}
\begin{pmatrix} Id_{\mathcal{H}_C}&0\\
		0&U_1\end{pmatrix},
\end{equation}
where $U_1$ is an unitary operator on $\mathcal{H}_1$.
if $g$ is represented by operator $U$ 
then
\begin{equation}\label{diag1}
U\begin{pmatrix} Id_{\mathcal{H}_C}&0&0\\
		0&U_1^{-1}&0\\
0&0&Id_{\mathcal{H}_2}\end{pmatrix}=
\begin{pmatrix} Id_{\mathcal{H}_C}&0&0\\
		0&Id_{\mathcal{H}_1}&0\\
0&0&U_2\end{pmatrix}
,
\end{equation}
where $U_2$ is an unitary operator on $\mathcal{H}_2$.
We have
\begin{equation}\label{diag1}
U=
\begin{pmatrix} Id_{\mathcal{H}_C}&0&0\\
		0&U_1&0\\
0&0&U_2\end{pmatrix}
.
\end{equation}

Operator $D_A|_{\mathcal{H}_2}$ commute with $U_2$.
Operator $D_A|_{\mathcal{H}_1}=D_B|_{\mathcal{H}_1}$ commute with $U_1$. 
Hence $D_A$ commute with $U$. We can say the same about $J_A$,
$\Gamma_A$.

Axiom 6.
Dimension of $\mathbf{A}$ equals to dimension of $\mathbf{B}$ and dimension of $\mathbf{B}$ equals to dimension of $\mathbf{C}$.
Hence dimension of $\mathbf{A}$ equals to dimension of $\mathbf{C}$.

Axiom 7.
Let $c_A$, $c_B$, $c_C$ are fundamental cycles on $\mathbf{A}$,  $\mathbf{B}$, $\mathbf{C}$.  
Then we have $c_A=f_1^*(c_B)$ and $c_B=g_1^*(c_C)$. Hence $c_A=(f_1g_1)^*(c_C)$. 
  
\end{proof}

\subsection{Fundamental group}
Let $\mathbf{A}$ be a spectral triple. Consider the following category. Object of this category is a finite covering $f_i:\mathbf{A_i}\rightarrow \mathbf{A}$.
Every  object of this category defines a group $G(\mathbf{A_i}, \mathbf{A})$. Morphism from $f_i:\mathbf{A_i}\rightarrow \mathbf{A}$ to $f_j:\mathbf{A_j}\rightarrow \mathbf{A}$
is a such finite covering $f_{ij}:\mathbf{A_i}\rightarrow \mathbf{A_j}$ that the diagram: 
\[
\begin{diagram}
	\node{\mathbf{A_i}} \arrow{se,t}{f_i} \arrow[2]{e,t}{f_{ij}} \node[2]{\mathbf{A_j}}\arrow{sw,t}{f_j} \\
	\node[2]{\mathbf{A}}
\end{diagram}
\]
is commutative.

Every morphism of this category naturally defines a surjective group homorphism 
$G(\mathbf{A_i}, \mathbf{A})\rightarrow G(\mathbf{A_j}, \mathbf{A})$.

Hence for every spectral triple we have a commutative diagram of groups and surjective homomorphisms.
\begin{defn}\label{defn-fund}
Fundamental group of spectral triple $\mathbf{A}=(\A, \mathcal{H}, D)$ is an inverse limit of described above diagram of groups.
\end{defn}
We shall use following notation $\pi_1(\mathbf{A})$ or $\pi_1(\A, \mathcal{H}, D)$ for fundamental group of $\mathbf{A}=(\A, \mathcal{H}, D)$.
\section{Fundamental group of commutative spectral triple}

Let us recall some facts of noncommutative and differential geometry and topology.
In \cite{varilly:noncom} and \cite{suprsym:qt} it was shown that every commutative spectral triple 
$\mathbf{B}=(\B, \mathcal{H}_B, D_B, J_B(,\Gamma_B))$ defines compact Riemann manifold $N$,
and spinor $\mathcal{S}_N$ bundle on it. In this case $B\approx C^\infty(N)$ and $\mathcal{H}_B\approx L^2(\mathcal{S}_N)$. Moreover, 
 $D_B, J_B(,\Gamma_B)$
correspond to local ``operators" on  smooth sections of the bundle $\mathcal{S}_N$. Accorging to \cite{varilly:noncom} the fundametal cycle on $\mathbf{B}$ corresponds to the volume form 
$\Omega_N$ of the Riemann manifold $N$. Inversely, if we have a Riemann manifold $N$ with spinor bundle $\mathcal{S}_N$ and local ``operators" on the bundle smooth sections
which satisfy to the set of conditions, then we can build corresponding spectral triple. 
A completion of a $pre-C^*$ algebra $\A$ is the $C^*$ algebra $A$, and latter defines
a topological space. 
Every *- homomorphism of commutative $C^*$ algebras defines 
continous map of topological spaces \cite{murphy}. 
In this case every element of $pre-C^*$ algebra corresponds to a smooth complex function. 
Hence if a *- homomorphism is an extention of a 
*-homomorphism of $pre-C^*$ algebras then the map is smooth. 
If $f:M\rightarrow N$ is a connected finitely sheeted covering of compact Riemann manifold $N$ then $M$ has the natural structure of 
Riemann manifold that is compact.
If a local diffeomorphism $f:M\rightarrow N$ is not surjective then $C^{\infty}(M)$ is not a finitely generated $C^{\infty}(N)$ module 
(we consider the natural $C^{\infty}(N)$ module structure). 

These facts will be used below to show that the fundametal group of commutative spectral triple $\mathbf{B}=(C^\infty(N), L^2(\mathcal{S}_N), D_B, J_B(,\Gamma_B))$
is the profinite completion of $\pi_1(N)$. 
\begin{lem}\label{comcov-first}
Let $\mathbf{B}=(\B, \mathcal{H}_B, D_B, J_B(,\Gamma_B))$ be a commutative spectal triple that corresponds to Riemann manifold $N$, and
$f:M\rightarrow N$ is finitely sheeted covering. Then there exist such natural commutative spectal triple

$\mathbf{A}=(\A, \mathcal{H}_A, D_A, J_A(,\Gamma_A))$
and finite covering $f^*:\mathbf{A}\rightarrow \mathbf{B}$ that $\A=C^\infty(M)$ and $f^*$ naturally corresponds to the smooth map $f$.  
\end{lem}

\begin{proof}
Let $\A=C^\infty(M)$ and $\mathcal{S}_M$ is a spinor bundle on $M$ that is the inverse image of $\mathcal{S}_N$.
Since $D_B$, $J_B$(,$\Gamma_B$) correspond to local ``operators" on  smooth sections of the bundle $\mathcal{S}_N$ we can naturally define 
similar ``operators" on  smooth sections of the bundle $\mathcal{S}_M$. Then we can define $\mathcal{H}_A=L^2(\mathcal{S}_M)$ and ``operators"
$D_A$, $J_A$(,$\Gamma_A$) on $\mathcal{H}_A$. A direct checking shows that  that $\mathbf{A}=(\A, \mathcal{H}_A, D_A, J_A (,\Gamma_A))$
satisfies to axioms of spectral triple. There are the natural injective *- homomorphism $f_1^*:\B\rightarrow \A$ and the natural map from the space of smooth sections of 
$\mathcal{S}_N$ to the space of smooth sections of $\mathcal{S}_M$. The later defines the injective homomorhism $f_2^*:\mathcal{H}_B\rightarrow\mathcal{H}_A$.
Every fixed on $f_1^*(\B)$ * - automorpsism of $\A$ corresponds to such homeomorpsism 
$\alpha :M\rightarrow M$ that following diagram 

\[
\begin{diagram}
	\node{M} \arrow{se,t}{f} \arrow[2]{e,t}{\alpha} \node[2]{M}\arrow{sw,t}{f} \\
	\node[2]{N}
\end{diagram}
\]
is commutative. Indeed $\alpha$ is isometry and number of similar isometries is finite since $f$
is a finitely listed covering. Every similar isometry generates fixed on image of $\mathcal{H}_B$ unitary
operator on $\mathcal{H}_A$. So we have the finite group $G(\mathbf{A}, \mathbf{B})$ that is naturally isomorphic to described above group of isometries. We have checked
condition (i) of \ref{defn-per}. The (ii) is evident, and (iii) contains axioms. Let us check
the axioms.

Axiom 1. Follows from \cite{karoubi:k} (Exercise II 6.6).

Axiom 2. According to \cite{karoubi:k} the space of sections of $\mathcal{S}_A$ is the tensor
product of $\A$ and space of sections of $\mathcal{S}_B$. Hence we have $\mathcal{H}_A \approx \A\otimes_{\B} \mathcal{H}_B$.

Axiom 3. Let us select a point $x_0\in M$. Let be such $U\subset M$ ``fundamental domain" that:

a) if $f_1(x_1)=f_1(x_2)$ and $x_1\in U$ then distance between $x_2$ and $x_0$ does not exceeds distance
between $x_1$ and $x_0$,

b) $U$ is the maximal subset of $M$ with such properties.

In this case $U$ is an open subset of $M$ and $f_1(U)$ is an open dense subset of $N$.  
Let $\alpha$ is isometry of $M$. We define isometry $\beta$ in the following way:

\[
\begin{diagram}
	\node{f_1(U)}\arrow{e,t}{\alpha|_U} \node{\alpha(U)} \arrow{e,t}{f_1} \node{N}  
\end{diagram}
.
\]

Since $f(U)$ is dense, the $\beta$ domain may be uniquely extended to $N$. So for
every isometry of $M$ we have the isometry of $N$. It gives the homorphism
$G(\mathbf{A})\rightarrow G(\mathbf{B})$.
Now we will show that the homomrphism is surjective. Let $\beta$ is isometry of $N$. Select such point $x'_0\in M$,
that $f(x'_0)=\beta (f(x_0)$. Select
pathwise connected, simple connected neighborhood $V$ of $f_1(x_0)$. Select such $U$ and $U'$ connected 
components of $f^{-1}(V)$ and $f^{-1}(\beta(V))$ that $x_0\in U$ and $x'_0\in U'$.
Since $f$ is a covering we can uniquely define such isometry $\alpha:U\rightarrow U'$ that
the following diagram
\[
\begin{diagram}
	\node{U} \arrow[2]{e,t}{\alpha} \arrow{s,l}{f}
								\node[2]{U'} \arrow{s,r}{f} \\
	\node{V} \arrow[2]{e,b}{\beta}	\node[2]{\beta(V)}
\end{diagram}
\]
is commutative. There are no obstacles to extend isometry $\alpha$ to $M$.
So every isometry of $N$ is the image of isometry of $M$, and homomorpsim $G(\mathbf{A})\rightarrow G(\mathbf{B})$
is surjective.

Axiom 4.

Follows from the construction of $G(\mathbf{A}, \mathbf{B})$.
  
Axiom 5.

Follows from the construction of $D_A$, $J_A$ (and $\Gamma_A$).

Axiom 6.

According to \cite{varilly:noncom}, dimension of commutative spectral
triple equals to dimension of correspond Riemann manifold. It is
clear that dimension of $M$ equals to dimension of $N$.

Axiom 7.

As it has been shown above the fundamental cycle on $\mathbf{A}$ ($\mathbf{B}$) corresponds to volume form
$\Omega_N$ ($\Omega_M$). Since metric on $M$ is the natural image of metric on $N$ the
$\Omega_M$  is natural image of $\Omega_N$ i. e. $\Omega_M=f^*\Omega_N$.
Then we have $c_A=f^*_1(c_B)$ 

\end{proof}

\begin{lem}\label{comcov-second}
Let $\mathbf{B}=(\B, \mathcal{H}_B, D_B, J_B(,\Gamma_B))$ be a commutative spectal triple 
and $f:\mathbf{A}= (\A, \mathcal{H}_A, D_A, J_A(,\Gamma_A))\rightarrow\mathbf{B}$ is a finite covering.
Then $\B$ belongs to the center of $\A$. 
\end{lem}

\begin{proof}
Suppose that $\B$ does not belong to the center of $\A$. As vector space $\B$ is
generated by elements $e^{i\psi}$ where $\psi$ is a real smooth function on the Riemann
manifold that corresponds to the $\mathbf{B}$. Hence we have one function $e^{i\psi}$ that does
not belong to the center of $\A$. For every $n\in\mathbf{N}$ we have unitary element $u_n=e^{i\psi/n}\in\B$
and cooresponding fixed on the $\B$ inner authomorphism of $\A$.
We have an infinite set of different inner automorphisms of $\A$ those are fixed on $\B$. However
every inner authomorphism corresponds to geometry equivalence \cite{varilly:noncom}.
Hence the group $G(\mathbf{A}, \mathbf{B})$ is infinite. It contradicts with \ref{defn-per}.

\end{proof}

\begin{lem}\label{comcov-third}
Let $\mathbf{B}=(\B, \mathcal{H}_B, D_B, J_B(,\Gamma_B))$ is a commutative spectal triple 
and 
$f:\mathbf{A}= (\A, \mathcal{H}_A, D_A, J_A(,\Gamma_A))\rightarrow\mathbf{B}$ is a finite covering.
Then $\A$ is commutative. 
\end{lem}

\begin{proof}
$\A$ is
generated by sefajoint elements as vector space.
Suppose that $\A$ is not commutative.  Then there existsts two selfajoint elements $x,y\in\A$ those do not commute.
It may be shown that among inner automorphisms that correspond
to unitary elements $e^{i\epsilon(x+ay)}$ ($a,\epsilon\in\mathbf{R}$)
there exists an infinite set of different ones.
Since $\B$ is contained in the center of $\A$ these automorphisms are fixed on $\B$.
Hence the group $G(\mathbf{A}, \mathbf{B})$ is infinite.
It contradicts with \ref{defn-per}.

\end{proof}

\begin{lem}\label{comcov-four}
Let $\mathbf{B}=(\B, \mathcal{H}_B, D_B, J_B(,\Gamma_B))$ be a commutative spectal triple 
that corresponds to Riemann manifold $N$, and $f:\mathbf{A}= (\A, \mathcal{H}_A, D_A, J_A(,\Gamma_A))\rightarrow\mathbf{B}$ is a finite covering.
Then

(i) There exists Riemann manifold $M$  and finitely listed covering 
$f:M\rightarrow N$ that
$\mathbf{A}=(C^{\infty}(M), L^2(f^*(\mathcal{S}_M)), D_A, J_A(,\Gamma_A))$ where $\mathcal{S}_M$ is inverse image of $\mathcal{S}_N$.

(ii) ``Operators" $D_A$, $J_A$ $(,\Gamma_A)$ are naturally defined by ``operators" $D_B$, $J_B$ $(,\Gamma_B)$
on smooth sections of $\mathcal{S}_N$. 
\end{lem}
\begin{proof}
(i) According to \ref{comcov-third} $\A$ is commutative. Hence

$\mathbf{A}=(C^{\infty}(M), L^2(\mathcal{S}_M), D_A, J_A(,\Gamma_A))$.
Homomorphism $f_1$ difines the smooth map $f:\B\rightarrow \A$.
The axiom 7 means that $\Omega_M=f^*\Omega_N$ where $\Omega_M$ and $\Omega_N$ are volume
forms on $M$ and $N$. Hence Jacobian of $f$ has no zeros and $f$ is a local diffeomorphism. Since $\A$ is finitely generated $\B$ module
the $f$ is surjective. Every surjective local diffeomorfism is a covering. According to axiom 2 of finite covering 
$\mathcal{H}_A \approx \A\otimes_{\B} \mathcal{H}_B$. It means that $\mathcal{S}_M$ is the inverse image of $\mathcal{S}_N$.
(ii) This property follows from the axiom 5 of finite covering.
\end{proof}

The lemmas \ref{comcov-first} and  \ref{comcov-four} mean that there is
 one to one correspondence between finite coverings of 
commutative spectral triple 

$\mathbf{A}=(C^{\infty}(M), L^2(\mathcal{S}_M), D_A, J_A(,\Gamma_A))$
and finitely sheeted coverings of $M$.
For every finitely listed covering $f_i:M_i\rightarrow M$
let us call $G(M_i|M)$ the covering group \cite{spanier:at}, i. e. the group of such homeomorphisms of $M_i$
$\alpha$ that $f_i\alpha=f_i$. Then corresponding groups $G(M_i|M)$ and $G(\mathbf{A_i}, \mathbf{A})$
are naturally isomorphic. 
Moreover, according 
to the following diagrams:

\[
\begin{diagram}
	\node{\mathbf{A_i}} \arrow{se,t}{f_i} \arrow[2]{e,t}{f_{ij}} \node[2]{\mathbf{A_j}}\arrow{sw,t}{f_j} \\
	\node[2]{\mathbf{A}}
\end{diagram}
\]
\[
\begin{diagram}
	\node{M_i} \arrow{se,t}{f^*_i} \arrow[2]{e,t}{f^*_{ij}} \node[2]{M_j}\arrow{sw,t}{f^*_j} \\
	\node[2]{M}
\end{diagram}
\]
there exists one to one correspondence between finite coverings 
$f_{ij}:\mathbf{A_i}\rightarrow\mathbf{A_j}$ and finitely generated coverings $f^*_{ij}:M_i\rightarrow M_j$.
Moreover this corresondence implyies that we have the following commutative diagram:
\[
\begin{diagram}
	\node{G(M|M_i)} \arrow[2]{e,t}{\approx} \arrow{s,l}{}
								\node[2]{G(\mathbf{A_i}, \mathbf{A})} \arrow{s,r}{} \\
	\node{G(M|M_j)} \arrow[2]{e,b}{\approx}	\node[2]{G(\mathbf{A_j}, \mathbf{A})}
\end{diagram}
\]

Let $\wt{M}$ be the universal covering of $M$. According to \cite{spanier:at}
$\pi_1(M)$ acts on $\wt{M}$ and $M\approx\wt{M}/\pi_1(M)$.
If $f^*_i:M_i\rightarrow M$ is finitely sheeted covering then $M_i\approx\wt{M}/H_i$ and 
$G(M_i|M)\approx \pi_1(M)/H_i$ where $H_i$ is a finite index normal subroup of $\pi_1(M)$.
Inversely, every finite factorgroup difines a finitely sheeted covering.
Recall that profinite completion $\wh{G}$ of group $G$ is the inverse
limit of the diagram of its finite factorgroups \cite{johnstone:topos}.

\begin{thm}
If $\mathbf{A}=(C^{\infty}(M), L^2(\mathcal{S}_M), D_A, J_A(,\Gamma_A))$ is a commutative
spectral triple then $\pi_1 (\mathbf{A})\approx\wh{\pi_1(M)}$.
\end{thm}

\begin{proof}
Consider
\[
\begin{diagram}
	\node{G_i} \arrow[2]{e,t}{\alpha_{ij}} \node[2]{G_j}\\
\end{diagram}
\]
the diagram of finite factorgroups of $\pi_1(M)$. As it has been explained above
$G_i\approx G(M_i|M)$ for some finitely shifted covering.

We have the following diagram of finiely shifted coverings

\[
\begin{diagram}
	\node{M_i} \arrow{se,t}{f^*_i} \arrow[2]{e,t}{f^*_{ij}} \node[2]{M_j}\arrow{sw,t}{f^*_j} \\
	\node[2]{M}
\end{diagram}
\]

There is one to one correspondence of objects and arrows of this diagram to following diagram

\[
\begin{diagram}
	\node{\mathbf{A_i}} \arrow{se,t}{f_i} \arrow[2]{e,t}{f_{ij}} \node[2]{\mathbf{A_j}}\arrow{sw,t}{f_j} \\
	\node[2]{\mathbf{A}}
\end{diagram}
\]

In fact, both the diagram of spectral triples and the diagram of Riemann manifolds correspond to the same
diagram of groups and surjective homomorphisms.

According to definition \ref{defn-fund} $\pi_1(\mathbf{A})$ is the inverse limit of the diagram
\[
\begin{diagram}
	\node{G(\mathbf{A_i}, \mathbf{A})} \arrow[2]{e,t}{} \node[2]{G(\mathbf{A_j}, \mathbf{A})}\\
\end{diagram}
\]
Hence, it is
isomorphic to the diagram of finite factorgroups of $\pi_1(M)$. Thus
$\pi_1 (\mathbf{A})$ is the profinite completion of $\pi_1(M)$.

\end{proof}

\section{Aknowlegment} 
Authors would like to aknowldge members of Russian Federation Seminar ``Noncommutative Geometry and Topology"
leaded by professsor A. S. Mischenko and others
for discussion of this work.

\end{document}